\documentclass{elsart}

\usepackage{graphics}
\usepackage{graphicx}

\usepackage{amssymb}
\usepackage{amsmath}

\begin{document}

\begin{frontmatter}



\title{On existence and uniqueness of solution for a hydrodynamic problem 
related to water artificial circulation in a lake\thanksref{label1}
}

\thanks[label1]{Supported by project MTM2015-65570-P of MINECO/FEDER (Spain).}


\author[Santiago]{Francisco J. Fern\'andez},
\ead{fjavier.fernandez@usc.es}
\author[Vigo]{Lino J. Alvarez-V\'azquez\corauthref{cor}},
\corauth[cor]{Corresponding author. Tel: +34 986 812166. Fax: +34 986 812116.}
\ead{lino@dma.uvigo.es}
\author[Vigo]{Aurea Mart\'{\i}nez}
\ead{aurea@dma.uvigo.es}

\address[Santiago]{Universidade de Santiago de Compostela, Fac. Matem\'aticas, 15782 Santiago, Spain.}
\address[Vigo]{Universidade de Vigo, E.I. Telecomunicaci\'on, 36310 Vigo, Spain.}

\begin{abstract}
In this work we introduce a well-posed mathematical model for the processes involved in the artificial circulation of water, in order to avoid eutrophication phenomena, for instance, in a lake. 
This novel and general formulation is based on the modified Navier-Stokes equations following the Smagorinsky model of turbulence, and presenting a suitable nonhomogeneous Dirichlet boundary condition.
For the analytical study of the problem, we prove several theoretical results related to existence, uniqueness and smoothness for the solution of this recirculation model.
\end{abstract}

\end{frontmatter}


\section{Introduction}

Artificial circulation is an usual and effective technique for remediating eutrophic bodies of water suffering from oxygen depletion or algal blooms. 
By means of the mechanical mixing of water, natural stratification is broken, so that the creation of a water circulation patterns allows water aeration and,
consequently, an increasing in the dissolved oxygen content and a decreasing in water quality problems. 

The use of artificial circulation as a water aeration technique is based on forcing water to expose to the atmosphere (main source of oxygen through diffusion processes).
In order to create a circulation pattern preventing stratification, flow pumps take water from the upper layer of water (the well aerated epilimnion) by means of a collector, 
injecting it into the bottom layer (the poorly oxygenated hypolimnion). So, water from the bottom is circulated to the surface, where oxygenation from the atmosphere can naturally occur. 
Despite its obvious practical interest, as far as we know, the environmental problem has not been addressed before from the viewpoint of 
modelling and its mathematical analysis. The first steps in this direction, though from an optimal control perspective, have been recently given by the authors 
in the work \cite{mcrf}.

In the first part of this paper we focus our attention on setting a novel, general, well-posed mathematical model for water artificial circulation, based on the 
modified Navier-Stokes equations with a Smagorinsky term of turbulence, and completed with nonhomogeneous Dirichlet boundary conditions 
(translating the complex mechanical behaviour of the process). 
Although in the recent mathematical literature it is easy to find numerous articles devoted to the numerical treatment of the Smagorinsky model 
(without trying to be exhaustive, and only in the last decade, we can mention among others \cite{n0,n1,n2,n3,n4,n5,n6,n7,n8}), its analytical study has received much less attention. 
We must cite here the pioneering works \cite{lady1,Gunzburger2,t1,t2},
mainly devoted to the homogeneous boundary case, unfortunately not applicable here. 
Thus, in the second part of the paper we address a rigorous derivation of several analytical results aimed to the well-posedness of our model. 
So, in particular, we present the detailed proofs of existence, uniqueness and regularity for the solution of this recirculation model.

\section{Setting of the governing equations}

To fix ideas, we consider a domain $\Omega \subset \mathbb{R}^3$ corresponding, for instance, to a lake. 
In order to promote the artificial circulation of water inside the domain, we suppose the existence of a set of $N_{CT}$ pairs 
collector-injector in such a way that each water collector is connected to its corresponding 
injector by a pipe with a pumping group. We assume a smooth enough boundary $\partial \Omega$, such that it 
can be split into three disjoint subsets
$\partial \Omega = \Gamma_C \cup \Gamma_T \cup \Gamma_N$,
where $\displaystyle \Gamma_C$ corresponds to the part of the boundary where the water collectors $C^k, \, k=1,\ldots,N_{CT},$ are located,
$\displaystyle \Gamma_T$ corresponds to the part of the boundary where the water injectors $T^k, \, k=1,\ldots,N_{CT},$ are located, and
$\displaystyle \Gamma_N=\partial \Omega \setminus \left( \Gamma_C \cup \Gamma_T \right)$ corresponds to the rest of the boundary.

\begin{figure}[tb]
\begin{center}
\includegraphics[scale=0.95]{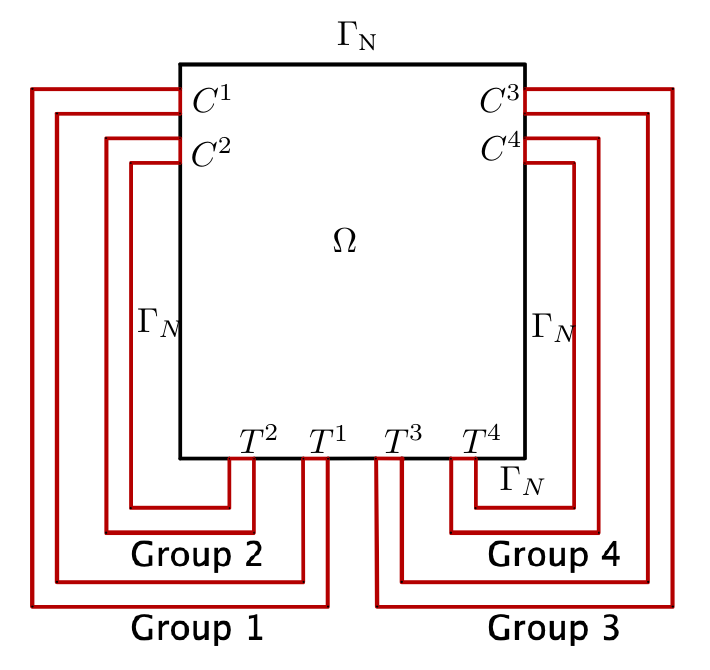}
\caption{Particular case of domain $\Omega$ with $N_{CT}=4$ injector-collector pairs.} \label{fig2}
\end{center}
\end{figure}

In order to simulate hydrodynamics we consider the modified Navier-Stokes equations following the Smagorinsky model of turbulence:
\begin{equation}\label{eq:system1}
\left\{\begin{array}{l}
\displaystyle
\frac{\partial \mathbf{v}}{\partial t} + \nabla \mathbf{v} \mathbf{v} - \nabla \cdot \Xi(\mathbf{v})
+\nabla p = \mathbf{F}
\quad \mbox{in} \; \Omega \times ]0,T[, \\
\displaystyle
\nabla \cdot \mathbf{v}=0 \quad \mbox{in} \; 
\Omega \times ]0,T[, \\
\displaystyle
\mathbf{v}=\boldsymbol{\phi}_{\mathbf{g}} \quad \mbox{on}\; 
\partial \Omega\times ]0,T[, \\
\displaystyle
\mathbf{v}(0)=\mathbf{v}_0 \quad \mbox{in}\; \Omega,
\end{array}\right.
\end{equation}
where $\mathbf{v}(\mathbf{x},t)$ is the velocity of water, $\mathbf{F}(\mathbf{x},t)$ stands for the source term, 
$\mathbf{n}$ represents the unit outward normal vector to the boundary 
$\partial \Omega$, $\mathbf{v}_0(\mathbf{x})$ is the initial velocity, and boundary function $\boldsymbol{\phi}_{\mathbf{g}}(\mathbf{x},t)$ 
is given by the following expression:
\begin{equation} \label{eq:phig}
\boldsymbol{\phi}_{\mathbf{g}}(\mathbf{x},t)
=\sum_{k=1}^{N_{CT}}g^k(t) \left[\frac{\varphi^k(\mathbf{x})}{\mu(T^k)}
 -\frac{\widetilde{\varphi}^k(\mathbf{x})}{\mu(C^k)} 
\right] \mathbf{n}.
\end{equation}
where $\mu(S)$ denotes the usual measure of a generic set $S$, for each $k=1,\ldots,N_{CT}$, 
$g^k(t) \in H^1(0,T)$ is the positive function representing the volumetric flow rate by pump $k$ at each time $t$, and mappings
$\varphi^k,\, \widetilde{\varphi}^k \in H^{3/2}(\partial \Omega)$, 
$k=1,\ldots,N_{CT}$, corresponding to collectors and injectors, satisfy the following assumptions:
\begin{itemize}
\item $\varphi^k(\mathbf{x}),\, \widetilde{\varphi}^k(\mathbf{x}) \geq 0 \quad  $
a.e. $\mathbf{x} \in \partial \Omega$, 
\item $\varphi^k(\mathbf{x})=0 \quad $ a.e. $\mathbf{x} \in \partial \Omega 
\setminus T^k$, and $\displaystyle \int_{T^k} \varphi^k(\mathbf{x}) \, d \gamma =\mu(T^k)$,
\item $\widetilde{\varphi}^k(\mathbf{x})=0 \quad$ a.e. $\mathbf{x} \in \partial \Omega 
\setminus C^k$, and $\displaystyle \int_{C^k} \widetilde{\varphi}^k(\mathbf{x})\, d \gamma = 
\mu(C^k)$.
\end{itemize}

\begin{rem}
It is clear that, under these assumptions, the boundary function $\boldsymbol{\phi}_{\mathbf{g}}$ verifies the usual compatibility condition 
$\int_{\partial \Omega} \boldsymbol{\phi}_{\mathbf{g}}(\mathbf{x},t) \; d \gamma = \boldsymbol{0}$. 
One of the simplest examples of functions $\varphi^k$ and $ \widetilde{\varphi}^k$, $k=1,\ldots,N_{CT},$
satisfying above hypotheses is given by a suitable regularization of the indicator functions $1_{T^k}$ and $ 1_{C^k}$, respectively.
\end{rem}

Finally, the turbulence term $\Xi(\mathbf{v})$ is given by:
\begin{equation}
\Xi(\mathbf{v})=\left. \frac{\partial D(\epsilon)}{\partial \epsilon} \right|_{\epsilon=\epsilon(\mathbf{v})},
\ \mbox{ with } \epsilon(\mathbf{v})=\frac{1}{2} \left( \nabla \mathbf{v} + \nabla \mathbf{v}^t\right),
\end{equation}
where $D$ is a potential function (for instance, in the particular case of the classical Navier-Stokes equations,
$D(\epsilon)=\nu \left[\epsilon:\epsilon\right]$ and, consequently,
$\Xi(\mathbf{v}) = 2 \nu\, \epsilon(\mathbf{v})$).
In our case, the Smagorinsky turbulence model, the potential function $D$ is defined as \cite{lady1}:
\begin{equation} \label{eq:doperador}
D(\epsilon)=\nu \left[\epsilon:\epsilon\right] + \frac{2}{3} \nu_{tur} \left[\epsilon:\epsilon\right]^{3/2}.
\end{equation}
So,
\begin{equation}
\begin{array}{rcl}
\displaystyle
\Xi(\mathbf{v}) &=&
\displaystyle \left. \frac{\partial D(\epsilon)}{\partial \epsilon} \right|_{\epsilon=\epsilon(\mathbf{v})} = \displaystyle
2 \nu \, \epsilon(\mathbf{v}) + 2 \nu_{tur}  \left[\epsilon(\mathbf{v}):\epsilon(\mathbf{v})\right]^{1/2}
\epsilon(\mathbf{v}) \\
&=& \displaystyle
\left(2 \nu+ 2 \nu_{tur}  \left[\epsilon(\mathbf{v}):\epsilon(\mathbf{v})\right]^{1/2} \right) \epsilon(\mathbf{v}) = \displaystyle
\beta(\epsilon(\mathbf{v})) \, \epsilon(\mathbf{v}),
\end{array}
\end{equation}
where $\beta(\epsilon(\mathbf{v}))=
2 \nu+ 2 \nu_{tur}  \left[\epsilon(\mathbf{v}):\epsilon(\mathbf{v})\right]^{1/2} $. 

From a mathematical point of view, the advantage of considering the modified 
Navier-Stokes equations, besides being more appropriate for turbulent flows, 
lies in the fact that, if the potential function fulfills certain properties (see, for 
instance, related comments in \cite{Gunzburger1}), it is 
possible to demonstrate the uniqueness of solution, in addition to gain in 
regularity with respect to the smoothness obtained for the original Navier-Stokes equations.
The potential function (\ref{eq:doperador}) for the Smagorinsky turbulence model 
considered in this work, as we will note when analyzing the existence of solution for the 
state equations, fulfills the necessary conditions to guarantee the uniqueness and 
regularity of the solution.

\section{The concept of solution}

We start this section by defining the functional spaces used in the search for solution 
of the system (\ref{eq:system1}). So, we consider:
\begin{equation}
\begin{array}{rcl}
\displaystyle
\mathbf{X}&=&
\displaystyle
\left\{\mathbf{v}\in [W^{1,3}(\Omega)]^3:\;
\nabla \cdot \mathbf{v}=0 , \quad  \mathbf{v}_{|_{\Gamma_N}}= \boldsymbol{0}
\right\} , \\
\displaystyle
\widetilde{\mathbf{X}}&=& \displaystyle \left\{\mathbf{v}\in [W^{1,3}(\Omega)]^3:\;
\nabla \cdot \mathbf{v}=0, \quad \mathbf{v}_{|_{\partial \Omega}}= \boldsymbol{0}
\right\} .
\end{array}
\end{equation}

In order to define an appropriate space for the solution of problem (\ref{eq:system1}), 
we consider, for a Banach space $V_1$ and a locally convex space $V_2$ such that $V_1\subset V_2$, 
the following Sobolev-Bochner space (cf. Chapter 7 of \cite{Roubicek1}), for $1\leq p,q \leq \infty$:
\begin{equation}
W^{1,p,q}(0,T;V_1,V_2)=\left\{u \in L^p(0,T;V_1):\; \frac{du}{dt} \in L^q(0,T;V_2)\right\},
\end{equation}
where $\frac{du}{dt} $ denotes the derivative of $u$ in the sense of distributions. It is well known that, if both $V_1$ and $V_2$ are Banach 
spaces, then $W^{1,p,q}(0,T;V_1,V_2)$ is also a Banach space endowed 
with the norm $\|u\|_{W^{1,p,q}(0,T;V_1,V_2)}=\|u\|_{L^p(0,T;V_1)}+
\left\| \frac{du}{dt}  
\right\|_{L^{q}(0,T;V_2)}$. 

Then, we define the following spaces that will be used in 
the mathematical analysis of system (\ref{eq:system1}):
\begin{equation}
\begin{array}{rcl}
\displaystyle
\mathbf{W}&=&
\displaystyle
W^{1,\infty,2}(0,T;\mathbf{X},[L^2(\Omega)]^3) \cap 
\mathcal{C}([0,T];\mathbf{X}), \\
\displaystyle
\widetilde{\mathbf{W}}&=&
\displaystyle
W^{1,\infty,2}(0,T;\widetilde{\mathbf{X}},[L^2(\Omega)]^3) \cap 
\mathcal{C}([0,T];\widetilde{\mathbf{X}}).
\end{array}
\end{equation}
We will assume that source term $\mathbf{F} \in L^2(0,T;[L^2(\Omega)]^3)$, initial condition $\mathbf{v}_0 \in \{\mathbf{v} \in 
W^{2,2}(\Omega):\; \nabla \cdot \mathbf{v}=0\; \mbox{and}\; \mathbf{v}_{|_{\partial \Omega}}= \boldsymbol{0}\}\subset 
\widetilde{\mathbf{X}}$ (this additional regularity will be necessary in the proof of the existence of solution), 
and function $\mathbf{g} = (g^1, \ldots , g^{N_{CT}})$ satisfies the compatibility condition $\mathbf{g}(0)= \boldsymbol{0}$. 
Now, we will prove an extension theorem that will be useful in the definition of solution. 

\begin{lem} \label{traza} There exists a linear continuous extension:
\begin{equation} 
\begin{array}{rcl}
R_{\mathbf{v}}: [H^1(0,T)]^{N_{CT}}& \rightarrow &
W^{1,2,2}(0,T;[H_{\sigma}^2(\Omega)]^3,[H_{\sigma}^2(\Omega)]^3) \\
\mathbf{g} & \rightarrow &R_{\mathbf{v}}(\mathbf{g} )= \boldsymbol{\zeta}_{\mathbf{g}},
\end{array}	
\end{equation}
such that ${\boldsymbol{\zeta}_{\mathbf{g}}}_{|_{\partial \Omega}}=
\boldsymbol{\phi}_{\mathbf{g}}$, where $\boldsymbol{\phi}_{\mathbf{g}}$ 
is defined by (\ref{eq:phig}), and $H_{\sigma}^2(\Omega)=\{ \mathbf{u} 
\in [H^2(\Omega)]^3:\; \nabla \cdot \mathbf{u}=0\}$.
\end{lem}

\noindent {\bf Proof } We denote, for $k=1,\ldots,N_{CT}$,
\begin{equation}
\boldsymbol{\psi}^k(\mathbf{x})=\left[\frac{\varphi^k(\mathbf{x})}{\mu(T^k)}
 -\frac{\widetilde{\varphi}^k(\mathbf{x})}{\mu(C^k)} 
\right] \mathbf{n} \in [H^{3/2}(\partial \Omega)]^3,
\end{equation}
we have that 
\begin{equation}
\int_{\partial \Omega} \boldsymbol{\psi}^k(\mathbf{x}) \cdot \mathbf{n} \, 
d \gamma = \frac{1}{\mu(T^k)}\int_{\partial \Omega}  \varphi^k(\mathbf{x}) \, d
\gamma- \frac{1}{\mu(C^k)} \int_{\partial \Omega} \widetilde{\varphi}^k(
\mathbf{x}) \, d \gamma=0.
\end{equation}
Thanks to results proved in \cite{Cattabriga}, there exists, for 
$k=1,\ldots,N_{CT}$, $\boldsymbol{\zeta}^k \in [H^2(\Omega)]^3$ 
and $p^k \in H^1(\Omega)/\mathbb{R}$ solutions of the following Stokes system:
\begin{equation}
\left\{\begin{array}{l}
\displaystyle
-\nu \Delta \boldsymbol{\zeta}^k+\nabla p^k=0 \quad \mbox{in}\; \Omega,  \\
\nabla \cdot \boldsymbol{\zeta}^k=0 \quad \mbox{in}\; \Omega,\\
\boldsymbol{\zeta}^k= \boldsymbol{\psi}^k  \quad \mbox{on}\; \partial \Omega,
\end{array}\right.
\end{equation}
that satisfy the following inequality:
\begin{equation}
\|\boldsymbol{\zeta}^k\|_{[H^2(\Omega)]^3}+\|p^k\|_{H^1(\Omega)/\mathbb{R}} \leq C
 \|\boldsymbol{\psi}^k \|_{[H^{3/2}(\partial \Omega)]^3}.
\end{equation}
Then, if we define
\begin{equation}
\boldsymbol{\zeta}_{\mathbf{g}}(t,\mathbf{x})=\sum_{k=1}^{N_{CT}} g^k(t) \boldsymbol{\zeta}^k(\mathbf{x}),
\end{equation}
we have that $\boldsymbol{\zeta}_{\mathbf{g}}\in W^{1,2,2}(0,T;[H^2(\Omega)]^3,[H^2(\Omega)]^3)$, and:
\begin{eqnarray*} \hspace*{-1.cm}
\left\|\boldsymbol{\zeta}_{\mathbf{g}} \right\|_{ W^{1,2,2}(0,T;[H^2(\Omega)]^3,[H^2(\Omega)]^3)} 
\leq C\sum_{k=1}^{N_{CT}}\|g^k\|_{H^1(0,T)} \left(
\frac{\|\varphi^k\|_{H^{3/2}(\partial \Omega)}}{\mu(T^k)}+
\frac{\|\widetilde{\varphi}^k\|_{H^{3/2}(\partial \Omega)}}{\mu(C^k)} \right),
\end{eqnarray*}
which concludes the proof. \hfill $\blacksquare$

\begin{lem} The following inclusion is continuous:
\begin{eqnarray}
&& W^{1,2,2}(0,T;[H^2(\Omega)]^3,[H^2(\Omega)]^3) \\
&& \qquad \subset 
W^{1,\infty,2}(0,T;[H^2(\Omega)]^3,[H^2(\Omega)]^3) \nonumber
\cap \mathcal{C}([0,T];[H^2(\Omega)]^3).
\end{eqnarray}
\end{lem}

\noindent {\bf Proof } We can apply Lemma 7.1 of \cite{Roubicek1}, taking 
$V_1=V_2=H^2(\Omega)$, and $p=q=2$.
\hfill $\blacksquare$

We consider now, for  $\mathbf{g} \in [H^1(0,T)]^{N_{CT}}$, the 
extension obtained in Lemma \ref{traza}:
\begin{equation}
\boldsymbol{\zeta}_{\mathbf{g}}\in 
W^{1,2,2}(0,T;[H^2(\Omega)]^3,[H^2(\Omega)]^3).
\end{equation}
We formally have that $\mathbf{v} \in \mathbf{W}$ is a solution of 
(\ref{eq:system1}) if and only if $\mathbf{z}=\mathbf{v}-\boldsymbol{\zeta}_{\mathbf{g}} 
\in \widetilde{\mathbf{W}}$ is a solution of the following system with 
homogeneous Dirichlet boundary conditions: 
\begin{equation} \label{eq:system2}
\left\{
\begin{array}{l}
\displaystyle \frac{\partial \mathbf{z}}{\partial t} + \nabla (\boldsymbol{\zeta}_{\mathbf{g}} +
\mathbf{z}) \mathbf{z} +\nabla \mathbf{z} \boldsymbol{\zeta}_{\mathbf{g}} 
 \vspace{.1cm} \\ \displaystyle  -\nabla \cdot \left(
2 \nu \epsilon(\mathbf{z})+
2 \nu_{tur} \int_{\Omega} \left[\epsilon(\boldsymbol{\zeta}_{\mathbf{g}}+\mathbf{z}):
\epsilon(\boldsymbol{\zeta}_{\mathbf{g}}+\mathbf{z}) \right]^{1/2} \epsilon(\boldsymbol{\zeta}_{\mathbf{g}}+\mathbf{z}) 
\right) +\nabla p   \\ \displaystyle
= \mathbf{F}- \frac{\partial \boldsymbol{\zeta}_{\mathbf{g}}}{\partial t}
-\nabla \boldsymbol{\zeta}_{\mathbf{g}} \boldsymbol{\zeta}_{\mathbf{g}} 
+2 \nu \nabla \cdot \epsilon(\boldsymbol{\zeta}_{\mathbf{g}})\quad \mbox{in} \; 
\Omega \times ]0,T[,  \\ \displaystyle
\mathbf{z}  = \mathbf{0} \quad \mbox{on} \; 
\partial \Omega \times ]0,T[,  \\
\displaystyle
\mathbf{z}(0)=\mathbf{v}_0 \quad \mbox{in}\; \Omega.
\end{array}
\right.
\end{equation}
Thus, we will use system (\ref{eq:system2}) to define the concept of solution 
for the system (\ref{eq:system1}). 

\begin{rem} It is worthwhile emphasizing here that:
\begin{itemize}
\item The sum $\boldsymbol{\zeta}_{\mathbf{g}} + \mathbf{z}$ is well defined in $\mathbf{W}$ since
$W^{1,2,2}(0,T;H^2(\Omega),H^2(\Omega))\subset \mathbf{W}$.
\item The term $2 \nu \nabla \cdot \epsilon(\boldsymbol{\zeta}_{\mathbf{g}})$ vanishes in the variational formulation 
thanks to the construction made in the proof of lemma \ref{traza}:
\begin{equation}
\nu \int_{\Omega} \epsilon( \boldsymbol{\zeta}_{\mathbf{g}} ): \epsilon( \boldsymbol{\eta})\, d \mathbf{x}=0,\quad
\forall  \boldsymbol{\eta} \in \widetilde{\mathbf{X}}.
\end{equation}
\item We must note that the regularity obtained in lemma \ref{traza} is more restrictive than the one
needed to define the concept of solution. However, this additional regularity for the extension will be necessary in next section in order to 
guarantee that the time derivative of the solution lies in $L^2(0,T;[L^2(\Omega)]^3)$. 
\end{itemize}
\end{rem}

Then, we can give the following definition of solution:

\begin{defn} \label{defsol} An element 
$\mathbf{v} \in \mathbf{W}$ is said to be a solution of problem (\ref{eq:system1}) if
there exists an element $\mathbf{z} \in \widetilde{\mathbf{W}}$ such that:
\begin{enumerate}
\item[a)] $\mathbf{v}=\boldsymbol{\zeta}_{\mathbf{g}} + \mathbf{z}$, with
$\displaystyle \boldsymbol{\zeta}_{\mathbf{g}} \in 
W^{1,2,2}(0,T;[H^2(\Omega)]^3,[H^2(\Omega)]^3)$
the reconstruction of the trace as given in Lemma \ref{traza}.
\item[b)] $\mathbf{z}(0)=\mathbf{v}_0 \ $ a.e. in $\Omega$.
\item[c)] $\mathbf{z}$ verifies the following variational formulation:
\begin{equation} \label{eq:system3}
\hspace{-1.55cm} \begin{array}{r}
\displaystyle
\int_{\Omega} \frac{\partial \mathbf{z}}{\partial t} \cdot \boldsymbol{\eta} \,d\mathbf{x} +
\int_{\Omega} \nabla (\boldsymbol{\zeta}_{\mathbf{g}}+
\mathbf{z}) \mathbf{z} \cdot \boldsymbol{\eta} \, d \mathbf{x} +
\int_{\Omega} \nabla \mathbf{z} \boldsymbol{\zeta}_{\mathbf{g}} 
\cdot \boldsymbol{\eta} \, d \mathbf{x}
+ 2 \nu \int_{\Omega} \epsilon(\mathbf{z}) : \epsilon(\boldsymbol{\eta}) \, 
d \mathbf{x} \\
+ \displaystyle
2 \nu_{tur} \int_{\Omega} \left[\epsilon(\boldsymbol{\zeta}_{\mathbf{g}}+\mathbf{z}):
\epsilon(\boldsymbol{\zeta}_{\mathbf{g}}+\mathbf{z}) \right]^{1/2} \epsilon(\boldsymbol{\zeta}_{\mathbf{g}}+\mathbf{z}):
\epsilon(\boldsymbol{\eta}) \, d \mathbf{x} \\
\displaystyle
=\int_{\Omega} \mathbf{H}_{\mathbf{g}} \cdot \boldsymbol{\eta} \, d \mathbf{x},
\quad \mbox{a.e.} \ t \in ]0,T[, \quad \forall
 \boldsymbol{\eta} \in \widetilde{\mathbf{X}},
\end{array}
\end{equation}
where
\begin{equation}
 \mathbf{H}_{\mathbf{g}}=\mathbf{F}-
\frac{\partial \boldsymbol{\zeta}_{\mathbf{g}}}{\partial t}
-\nabla \boldsymbol{\zeta}_{\mathbf{g}} \boldsymbol{\zeta}_{\mathbf{g}} \in L^2(0,T;[L^2(\Omega)]^3).
\end{equation}
\end{enumerate}
\end{defn}

\section{Existence of solution}

In order to better understand the results proved in this section, we will divide it in two parts: In the first subsection we 
will obtain some a priori estimates for the solution. In the second one, we will prove the existence of solution 
for the Galerkin approximations of system (\ref{eq:system2}), and we will demonstrate that the limit $\mathbf{v}$ of these 
Galerkin approximations is a solution of the system (\ref{eq:system2}) and that, consequently, the system (\ref{eq:system1})
admits, at least, a solution of the form $\mathbf{v}=\boldsymbol{\zeta}_{\mathbf{g}} + \mathbf{z}$.  
Uniqueness of solution will be analyzed in the last section.

\subsection{A priori estimates}

\begin{lem} If $\mathbf{z} \in \widetilde{\mathbf{W}}$ is a solution of the system (\ref{eq:system2}) in 
the sense of Definition \ref{defsol}, then there exist positive constants $C_1$ and $C_2$ such that:
\begin{equation} \label{eq:acotav1}
\begin{array}{r}
\displaystyle
\|\mathbf{z}\|^2_{L^{\infty}(0,T;[L^2(\Omega)]^3)}+
\|\mathbf{z}\|^2_{L^2(0,T;[W^{1,2}(\Omega)]^3)}+
\|\mathbf{z}\|^3_{L^3(0,T;[W^{1,3}(\Omega)]^3)}
\vspace{0.2cm} \\ \displaystyle
 \leq C_1  \bigg(
\|\mathbf{v}_0\|_{[L^2(\Omega)]^3}^2
+
\|\boldsymbol{\zeta}_{\mathbf{g}}\|^3_{L^3(0,T;[W^{1,3}(\Omega)]^3)}
+
 \|\mathbf{H}_{\mathbf{g}}\|^2_{L^2(0,T;[L^2(\Omega)]^2)}
\bigg).
\end{array}
\end{equation}
\begin{equation} \label{eq:acotav2}
\begin{array}{r}
\displaystyle
\left\|\frac{\partial \mathbf{z}}{\partial t} \right\|^2_{L^2(0,T;[L^2(\Omega)]^3)} + 
\|\mathbf{z}\|^2_{L^{\infty}(0,T;[W^{1,2}(\Omega)]^3)}+
\|\mathbf{z}\|^3_{L^{\infty}(0,T;[W^{1,3}(\Omega)]^3)} 
\vspace{0.2cm} \\  \displaystyle
\leq 
C_2 \Big(\nu \|\epsilon(\mathbf{v}_0) \|^2_{[L^2(\Omega)]^{3 \times 3}} + 
\frac{2}{3} \nu_{tur} \|\epsilon(\mathbf{v}_0) \|^3_{[L^3(\Omega)]^{3 \times 3}} 
+
\|\widetilde{\mathbf{H}}_{\mathbf{g}}\|^2_{L^2(0,T;[L^2(\Omega)]^3)} \Big)
\vspace{0.2cm}
\\
\displaystyle
\times \Bigg[1+
\exp \Big(\|\mathbf{v}_0\|_{[L^2(\Omega)]^3}^2
+
\|\boldsymbol{\zeta}_{\mathbf{g}}\|^3_{L^3(0,T;[W^{1,3}(\Omega)]^3)}
+ \|\mathbf{H}_{\mathbf{g}}\|^2_{L^2(0,T;[L^2(\Omega)]^2)}
\vspace{0.2cm} \\  \displaystyle
+\left\| \frac{\partial \boldsymbol{\zeta}_{\mathbf{g}}}{\partial t}\right\|^2_{L^2(0,T;H^1(\Omega))}
+\left\| \frac{\partial \boldsymbol{\zeta}_{\mathbf{g}}}{\partial t}\right\|^3_{L^2(0,T;W^{1,3}(\Omega))} \Big)
\Bigg]
\end{array}
\end{equation}
\end{lem}

\noindent {\bf Proof } For the sake of simplicity, we will divide the proof into two parts: So, in the first part we will take 
$\boldsymbol{\eta}=\mathbf{z}(t)$ as a test function in the variational formulation (\ref{eq:system3}) and 
we will obtain the first estimate. Finally, in the second part, we will take 
$\boldsymbol{\eta}=\frac{\partial \mathbf{z}}{\partial t}$ as a test function in order to obtain an 
estimate for the time derivative. 

As above commented, we first take $\boldsymbol{\eta}=\mathbf{z}(t)$ as a test function in the variational formulation (\ref{eq:system3}) and, 
integrating over the time interval $[0,t]$, we obtain:
\begin{equation}
\begin{array}{r}
\displaystyle 
\frac{1}{2}\|\mathbf{z}(t)\|^2_{[L^2(\Omega)]^3} + 
\int_0^t \int_{\Omega} \nabla \boldsymbol{\zeta}_{\mathbf{g}} \mathbf{z} \cdot \mathbf{z} \, d 
\mathbf{x} \, ds + 2 \nu \int_0^t \|\epsilon(\mathbf{z})\|^2_{[L^2(\Omega)]^{3 \times 3}} \, ds \vspace{0.2cm} \\ 
\displaystyle
+ 2 \nu_{tur} \int_0^t \|\epsilon(\boldsymbol{\zeta}_{\mathbf{g}}+\mathbf{z})\|^3_{[L^3(\Omega)]^{3 \times 3}} 
=\frac{1}{2} \|\mathbf{v}_0\|^2_{[L^2(\Omega)]^3} + 
\int_0^t \int_{\Omega} \mathbf{H}_{\mathbf{g}} \cdot \mathbf{z}\, d \mathbf{x} \, ds \vspace{0.2cm} \\ 
\displaystyle
+ 2\nu_{tur} \int_0^t \int_{\Omega}  \left[
\epsilon(\boldsymbol{\zeta}_{\mathbf{g}}+\mathbf{z}):\epsilon(\boldsymbol{\zeta}_{\mathbf{g}}+\mathbf{z})
\right]^{1/2} \epsilon(\boldsymbol{\zeta}_{\mathbf{g}}+\mathbf{z}) : \epsilon(\boldsymbol{\zeta}_{\mathbf{g}}) 
\, d \mathbf{x} \, ds,
\end{array}
\end{equation}
where we have used that $\displaystyle \int_{\Omega} \nabla \mathbf{z} \boldsymbol{\zeta}_{\mathbf{g}} \cdot 
\mathbf{z}\, d \mathbf{x}=\int_{\Omega} \nabla \mathbf{z} \mathbf{z} \cdot \mathbf{z}\, d \mathbf{x}=0$, and where
$C$ is a positive constant that can depend of Korn inequality. 
Now, if we take into account the inequality 
$ab\leq (1/2)a^2\epsilon+
(1/2)b^2\epsilon^{-1}$, we obtain that
\begin{displaymath}
\begin{array}{rcl}
\displaystyle 
\int_0^t \int_{\Omega} \mathbf{H}_{\mathbf{g}} \cdot \mathbf{z}\, d \mathbf{x} \, ds & \leq & 
\displaystyle 
\frac{1}{2}\epsilon_1 \int_0^t \|\mathbf{z}(s)\|_{[L^2(\Omega)]^3}
+\frac{1}{2\epsilon_1} \int_0^t \|\mathbf{H}_{\mathbf{g}}
\|^2_{[L^2(\Omega)]^3} \vspace{0.2cm} \\ &\leq & \displaystyle 
C \epsilon_1 \int_0^t \|\epsilon(\mathbf{z})\|^2_{[L^2(\Omega)]^{3 \times 3}} \, ds +\frac{1}{2\epsilon_1} \int_0^t \|\mathbf{H}_{\mathbf{g}}
\|^2_{[L^2(\Omega)]^3},
\end{array}
\end{displaymath}
forall $\epsilon_1>0$. In the other hand, $ab=(a \epsilon^{2/3})(b \epsilon^{-2/3}) \leq 
(2/3) a^{3/2} \epsilon + (1/3) b^3 \epsilon^{-2}$, and then:
\begin{equation}
\begin{array}{r}
\displaystyle \int_0^t \int_{\Omega}  \left[
\epsilon(\boldsymbol{\zeta}_{\mathbf{g}}+\mathbf{z}):\epsilon(\boldsymbol{\zeta}_{\mathbf{g}}+\mathbf{z})
\right]^{1/2} \epsilon(\boldsymbol{\zeta}_{\mathbf{g}}+\mathbf{z}) : \epsilon(\boldsymbol{\zeta}_{\mathbf{g}}) 
\, d \mathbf{x} \, ds \vspace{0.2cm} \\ 
\displaystyle \leq 
\int_0^t \|\epsilon(\boldsymbol{\zeta}_{\mathbf{g}}+\mathbf{z})\|^2_{[L^3(\Omega)]^{3 \times 3}} 
\|\epsilon(\boldsymbol{\zeta}_{\mathbf{g}})\|_{[L^3(\Omega)]^{3\times 3}}\, ds \vspace{0.2cm} \\ 
\displaystyle \leq \frac{2}{3} \epsilon_2
\int_0^t \|\epsilon(\boldsymbol{\zeta}_{\mathbf{g}}+\mathbf{z})\|^3_{[L^3(\Omega)]^{3 \times 3}} \, ds 
+ \frac{1}{3 \epsilon_2^2}
\int_0^t \|\epsilon(\boldsymbol{\zeta}_{\mathbf{g}})\|^3_{[L^3(\Omega)]^{3\times 3}}\, ds,
\end{array}
\end{equation}
for all $\epsilon_2>0$. Finally:
\begin{equation} \hspace*{-.5cm}
\begin{array}{r}
\displaystyle
\int_0^t \int_{\Omega} \nabla \boldsymbol{\zeta}_{\mathbf{g}} \mathbf{z} \cdot \mathbf{z} \, d 
\mathbf{x} \, ds \leq  \displaystyle
\int_0^t \|\nabla \boldsymbol{\zeta}_{\mathbf{g}}\|_{[L^3(\Omega)]^{3 \times 3}} 
\|\mathbf{z}\|^2_{[L^3(\Omega)]^3}\, ds \vspace{0.2cm} \\ 
\leq\displaystyle 
\frac{C}{\epsilon_3^2} \int_0^t \|\epsilon(\boldsymbol{\zeta}_{\mathbf{g}})\|^3_{[L^{3}(\Omega)]^{3\times 3}}\, ds
+C \epsilon_3 \int_0^t \|\epsilon(\mathbf{z})\|^3_{[L^3(\Omega)]^{3\times 3}}\, ds \vspace{0.2cm} \\ 
\leq  \displaystyle 
C \left(\frac{1}{\epsilon_3}+\epsilon_3 \right)
\int_0^t \|\epsilon(\boldsymbol{\zeta}_{\mathbf{g}})\|^3_{[L^{3}(\Omega)]^{3\times 3}}\, ds
+C \epsilon_3 \int_0^t \|\epsilon(\mathbf{z}+
\boldsymbol{\zeta}_{\mathbf{g}})\|^3_{[L^3(\Omega)]^{3\times 3}}\, ds,
\end{array} 
\end{equation}
for all $\epsilon_3>0$, where $C$ is a positive constant. If we adjust the 
values of $\epsilon_k$, $k=1,2,3$, we have that:
\begin{equation}
\begin{array}{r}
\displaystyle
\|\mathbf{z}(t)\|^2_{[L^2(\Omega)]^3}  
+\int_0^t \|\epsilon( \mathbf{z}(s))\|^2_{[L^2(\Omega)]^{3 \times 3}}\, ds
+\int_0^t \|\epsilon(\boldsymbol{\zeta}_{\mathbf{g}}+
\mathbf{z})\|^3_{[L^3(\Omega)]^{3 \times 3}}
\vspace{0.2cm} \\ \displaystyle 
\leq C \Big(
\|\mathbf{v}_0\|^2_{[L^2(\Omega)]^3}
+\int_0^t \|\mathbf{H}_{\mathbf{g}}
\|^2_{[L^2(\Omega)]^3}
\int_0^t \|\epsilon(\boldsymbol{\zeta}_{\mathbf{g}})\|^3_{[L^3(\Omega)]^{3\times 3}}\, ds \Big).
\end{array}
\end{equation}
So, there exists a positive constant $C_1$ such that: 
\begin{equation} 
\begin{array}{r}
\displaystyle
\|\mathbf{z}\|^2_{L^{\infty}(0,T;[L^2(\Omega)]^3)}+
\|\mathbf{z}\|^2_{L^2(0,T;[W^{1,2}(\Omega)]^3)}+
\|\mathbf{z}\|^3_{L^3(0,T;[W^{1,3}(\Omega)]^3)}
\vspace{0.2cm} \\ \displaystyle
 \leq C_1  \bigg(
\|\mathbf{v}_0\|_{[L^2(\Omega)]^3}^2
+
\|\boldsymbol{\zeta}_{\mathbf{g}}\|^3_{L^3(0,T;[W^{1,3}(\Omega)]^3)}
+
 \|\mathbf{H}_{\mathbf{g}}\|^2_{L^2(0,T;[L^2(\Omega)]^2)}
\bigg),
\end{array}
\end{equation}
and first estimate (\ref{eq:acotav1}) is stated.

Now, in this second part of the proof, we will consider the following notation
\begin{equation}
\widetilde{\mathbf{H}}_{\mathbf{g}}=\mathbf{F}
-\frac{\partial \boldsymbol{\zeta}_{\mathbf{g}}}{\partial t} \in L^2(0,T;[L^2(\Omega)]^3)
\end{equation}
and we will take $\boldsymbol{\eta}=\frac{\partial \mathbf{z}}{\partial t}$ as a test function in the 
variational formulation (\ref{eq:system3}). So, we have that:
\begin{equation}
\begin{array}{r}
\displaystyle 
\left\| \frac{\partial \mathbf{z}}{\partial t}(t)\right\|^2_{[L^2(\Omega)]^3} + 
\frac{d}{dt} \int_{\Omega} D(\epsilon(\boldsymbol{\zeta}_{\mathbf{g}}(t)+\mathbf{z}(t)))\, d \mathbf{x} = 
\int_{\Omega} \widetilde{\mathbf{H}}_{\mathbf{g}}
\cdot \frac{\partial \mathbf{z}}{\partial t}(t) \, d \mathbf{x}  
\vspace{0.2cm} \\ 
\displaystyle 
-\int_{\Omega} \nabla (\boldsymbol{\zeta}_{\mathbf{g}}(t)+\mathbf{z}(t))  (\boldsymbol{\zeta}_{\mathbf{g}}(t)+\mathbf{z}(t)) \cdot 
\frac{\partial \mathbf{z}}{\partial t}(t) \, d \mathbf{x}  \vspace{0.2cm} \\ 
\displaystyle +
\int_{\Omega} \beta  (\boldsymbol{\zeta}_{\mathbf{g}}(t)+\mathbf{z}(t)) \epsilon (\boldsymbol{\zeta}_{\mathbf{g}}(t)+\mathbf{z}(t)) :
\epsilon\left(\frac{\partial \boldsymbol{\zeta}_{\mathbf{g}}}{\partial t}(t)\right)\, d \mathbf{x},
\end{array}
\end{equation}	
where the potential $D$ is defined by formula (\ref{eq:doperador}) and $\beta$ is given by $\beta(\epsilon(\mathbf{v}))=
2 \nu+ 2 \nu_{tur}  \left[\epsilon(\mathbf{v}):\epsilon(\mathbf{v})\right]^{1/2}$. Now, thanks to Young inequality 
$ab=(\sqrt{2}a)(b/\sqrt{2})\leq a^2+b^2/4$, we obtain:
\begin{equation}
\int_{\Omega} \widetilde{\mathbf{H}}_{\mathbf{g}}(t) \cdot \frac{\partial \mathbf{z}}{\partial t} (t) \, dx \leq 
\|\widetilde{\mathbf{H}}_{\mathbf{g}}(t)\|^2_{[L^2(\Omega)]^3} + 
\frac{1}{4} \left\| \frac{\partial \mathbf{z}}{\partial t}(t)\right\|^2_{[L^2(\Omega)]^3},
\end{equation}
and also:
\begin{equation}
\begin{array}{r}
\displaystyle
\int_{\Omega} \nabla (\boldsymbol{\zeta}_{\mathbf{g}}(t)+\mathbf{z}(t))  (\boldsymbol{\zeta}_{\mathbf{g}}(t)+\mathbf{z}(t)) \cdot 
\frac{\partial \mathbf{z}}{\partial t}(t) \, d \mathbf{x} 
\vspace{0.2cm} \\ \displaystyle
\leq 
\| \nabla (\boldsymbol{\zeta}_{\mathbf{g}}(t) 
+\mathbf{z}(t))  (\boldsymbol{\zeta}_{\mathbf{g}}(t)+\mathbf{z}(t)) \|^2_{[L^2(\Omega)]^3} 
+ \frac{1}{4} \left\| \frac{\partial \mathbf{z}}{\partial t}(t)\right\|^2_{[L^2(\Omega)]^3}.
\end{array}
\end{equation}
Then, integrating in the time interval $]0,t[$:
\begin{equation} \hspace*{-1.cm}
\begin{array}{r}
\displaystyle \frac{1}{2} \left\|\frac{\partial \mathbf{z}}{\partial t} \right\|^2_{L^2(0,t;[L^2(\Omega)]^3)} + 
\int_{\Omega} D(\epsilon(\boldsymbol{\zeta}_{\mathbf{g}}(t)+\mathbf{z}(t)))\, d \mathbf{x}- 
\int_{\Omega} D(\epsilon(\boldsymbol{\zeta}_{\mathbf{g}}(0)+\mathbf{z}(0)))\, d \mathbf{x} \\  
\displaystyle 
\leq \|\widetilde{\mathbf{H}}_{\mathbf{g}}\|^2_{L^2(0,t;[L^2(\Omega)]^3)} + 
\|\nabla (\boldsymbol{\zeta}_{\mathbf{g}}-\mathbf{z})(\boldsymbol{\zeta}_{\mathbf{g}}-\mathbf{z}) \|_{L^2(0,t;[L^2(\Omega)]^3)}^2 \\ 
+
\displaystyle 
\int_0^t \int_{\Omega} \beta  (\boldsymbol{\zeta}_{\mathbf{g}}(s)+\mathbf{z}(s)) 
\epsilon (\boldsymbol{\zeta}_{\mathbf{g}}(s)+\mathbf{z}(s)) :
\epsilon\left(\frac{\partial \boldsymbol{\zeta}_{\mathbf{g}}}{\partial t}(s)\right)\, d \mathbf{x}\, ds,
\end{array}
\end{equation}
and we obtain that:
\begin{equation}
\hspace{-2.cm}\begin{array}{r}
\displaystyle \frac{1}{2} \left\|\frac{\partial \mathbf{z}}{\partial t} \right\|^2_{L^2(0,t;[L^2(\Omega)]^3)} + 
\nu \|\epsilon(\boldsymbol{\zeta}_{\mathbf{g}}(t)+\mathbf{z}(t))\|^2_{[L^2(\Omega)]^{3 \times 3}}  + 
\frac{2}{3} \nu_{tur} \|\epsilon(\boldsymbol{\zeta}_{\mathbf{g}}(t)+\mathbf{z}(t))\|^3_{[L^3(\Omega)]^{3 \times 3}} 
\\  \displaystyle \leq 
\nu \|\epsilon(\mathbf{v}_0) \|^2_{[L^2(\Omega)]^{3 \times 3}} + 
\frac{2}{3} \nu_{tur} \|\epsilon(\mathbf{v}_0) \|^3_{[L^3(\Omega)]^{3 \times 3}} 
+ \|\widetilde{\mathbf{H}}_{\mathbf{g}}\|^2_{L^2(0,t;[L^2(\Omega)]^3)}   \\  
\displaystyle 
+\|\nabla (\boldsymbol{\zeta}_{\mathbf{g}}-\mathbf{z})(\boldsymbol{\zeta}_{\mathbf{g}}-\mathbf{z}) \|_{L^2(0,t;[L^2(\Omega)]^3)}^2
\\  \displaystyle 
+\int_0^t \int_{\Omega} \beta  (\boldsymbol{\zeta}_{\mathbf{g}}(s)+\mathbf{z}(s)) 
\epsilon (\boldsymbol{\zeta}_{\mathbf{g}}(s)+\mathbf{z}(s)) :
\epsilon\left(\frac{\partial \boldsymbol{\zeta}_{\mathbf{g}}}{\partial t}(s)\right)\, d \mathbf{x}\, ds.
\end{array}
\end{equation}
We will focus our attention into the last two terms. Thanks to Poincar\'e and Holder inequalities:
\begin{equation}
\begin{array}{r}
\displaystyle  
\int_0^t \|\nabla (\boldsymbol{\zeta}_{\mathbf{g}}(s)+
\mathbf{z}(s)) (\boldsymbol{\zeta}_{\mathbf{g}}(s)
+\mathbf{z}(s))\|_{[L^2(\Omega)]^3}^2 \, d s \\ 
\displaystyle
\leq C \int_0^t \|\boldsymbol{\zeta}_{\mathbf{g}}(s)+\mathbf{z}(s)\|_{
[L^6(\Omega)]^3}^2
\|\nabla (\boldsymbol{\zeta}_{\mathbf{g}}(s)+\mathbf{z}(s)) \|_{[L^3(\Omega)]^{3\times 3}}^2
\\
\displaystyle
\leq C
\int_0^t \| \nabla (\boldsymbol{\zeta}_{\mathbf{g}}(s)+\mathbf{z}(s))\|^2_{[L^3(\Omega)]^{3 \times 3}}
\\
\displaystyle
\times \bigg( \|\nabla(\boldsymbol{\zeta}_{\mathbf{g}}(s)+\mathbf{z}(s))\|^2_{[L^{2\times 2}(\Omega)]^3}
+\| \nabla (\boldsymbol{\zeta}_{\mathbf{g}}(s)+\mathbf{z}(s))\|^3_{[L^3(\Omega)]^{3 \times 3}}
\bigg) \, ds, 
\\
\displaystyle
\leq C
\int_0^t \| \epsilon (\boldsymbol{\zeta}_{\mathbf{g}}(s)+\mathbf{z}(s))\|^2_{[L^3(\Omega)]^{3 \times 3}}
\\
\displaystyle
\times \bigg( \|\epsilon(\boldsymbol{\zeta}_{\mathbf{g}}(s)+\mathbf{z}(s))\|^2_{[L^{2\times 2}(\Omega)]^3}
+\| \epsilon (\boldsymbol{\zeta}_{\mathbf{g}}(s)+\mathbf{z}(s))\|^3_{[L^3(\Omega)]^{3 \times 3}}
\bigg) \, ds.
\end{array}
\end{equation}
On the other hand:
\begin{equation} \hspace{-2.1cm}
\begin{array}{r}
\displaystyle \int_0^t \int_{\Omega} \beta  (\boldsymbol{\zeta}_{\mathbf{g}}(s)+\mathbf{z}(s)) 
\epsilon (\boldsymbol{\zeta}_{\mathbf{g}}(s)+\mathbf{z}(s)) :
\epsilon\left(\frac{\partial \boldsymbol{\zeta}_{\mathbf{g}}}{\partial t}(s)\right)\, d \mathbf{x}\, ds 
\vspace{0.2cm} \\  \displaystyle = 
\nu \int_0^t \int_{\Omega} \epsilon(\boldsymbol{\zeta}_{\mathbf{g}}(s)+\mathbf{z}(s)):
\epsilon\left(\frac{\partial \boldsymbol{\zeta}_{\mathbf{g}}}{\partial t}(s)\right)\, d \mathbf{x}\, ds \vspace{0.2cm} \\ 
\displaystyle
+\nu_{tur} \int_0^t \int_{\Omega} [\epsilon(\boldsymbol{\zeta}_{\mathbf{g}}(s)+\mathbf{z}(s)):
\epsilon(\boldsymbol{\zeta}_{\mathbf{g}}(s)+\mathbf{z}(s))]^{1/2} 
\epsilon(\boldsymbol{\zeta}_{\mathbf{g}}(s)+\mathbf{z}(s)): 
\epsilon\left(\frac{\partial \boldsymbol{\zeta}_{\mathbf{g}}}{\partial t}(s)\right)\, d \mathbf{x}\, ds.
\end{array} \nonumber
\end{equation}
Thus, a simple computation shows us that:
\begin{equation}
\hspace{-1.cm}\begin{array}{r}
\displaystyle 
\nu \int_0^t \int_{\Omega} \epsilon(\boldsymbol{\zeta}_{\mathbf{g}}(s)+\mathbf{z}(s)):
\epsilon\left(\frac{\partial \boldsymbol{\zeta}_{\mathbf{g}}}{\partial t}(s)\right)\, d \mathbf{x}\, ds \vspace{0.2cm} \\  
\displaystyle \leq 
\nu \int_0^t \|\epsilon(\boldsymbol{\zeta}_{\mathbf{g}}(s)+\mathbf{z}(s))\|_{[L^2(\Omega)]^{3 \times 3}} 
\left\|
\epsilon\left(\frac{\partial \boldsymbol{\zeta}_{\mathbf{g}}}{\partial t}(s)\right)
\right\|_{[L^2(\Omega)]^{3 \times 3}} \, ds \vspace{0.2cm} \\  \leq 
\displaystyle \frac{\nu t}{4} + \int_0^t 
\|\epsilon(\boldsymbol{\zeta}_{\mathbf{g}}(s)+\mathbf{z}(s))\|_{[L^2(\Omega)]^{3 \times 3}}^2
\left\|
\epsilon\left(\frac{\partial \boldsymbol{\zeta}_{\mathbf{g}}}{\partial t}(s)\right)
\right\|_{[L^2(\Omega)]^{3 \times 3}}^2 \, ds
\leq 
\displaystyle \frac{\nu t }{4}  \vspace{0.2cm} \\  + \displaystyle \int_0^t 
\bigg(\|\epsilon(\boldsymbol{\zeta}_{\mathbf{g}}(s)+\mathbf{z}(s))\|_{[L^2(\Omega)]^{3 \times 3}}^2 
+\|\epsilon(\boldsymbol{\zeta}_{\mathbf{g}}(s)+\mathbf{z}(s))\|_{[L^3(\Omega)]^{3 \times 3}}^3 
\bigg)
\left\|
\epsilon\left(\frac{\partial \boldsymbol{\zeta}_{\mathbf{g}}}{\partial t}(s)\right)
\right\|_{[L^2(\Omega)]^{3 \times 3}}^2 \, ds,
\end{array} \nonumber
\end{equation}
and:
\begin{equation}
\hspace{-1.cm}\begin{array}{r}
\displaystyle 
\nu_{tur} \int_0^t \int_{\Omega} [\epsilon(\boldsymbol{\zeta}_{\mathbf{g}}(s)+\mathbf{z}(s)):
\epsilon(\boldsymbol{\zeta}_{\mathbf{g}}(s)+\mathbf{z}(s))]^{1/2} 
\epsilon(\boldsymbol{\zeta}_{\mathbf{g}}(s)+\mathbf{z}(s)): 
\epsilon\left(\frac{\partial \boldsymbol{\zeta}_{\mathbf{g}}}{\partial t}(s)\right)\, d \mathbf{x}\, ds \vspace{0.2cm} \\ 
\displaystyle \leq 
\nu_{tur}\int_0^t \|\epsilon(\boldsymbol{\zeta}_{\mathbf{g}}(s)+\mathbf{z}(s)) \|^2_{[L^3(\Omega)]^{3 \times 3}} 
\left\|
\epsilon\left(\frac{\partial \boldsymbol{\zeta}_{\mathbf{g}}}{\partial t}(s)\right)
\right\|_{[L^3(\Omega)]^{3 \times 3}} \, ds \vspace{0.2cm} \\ 
\displaystyle \leq \frac{4 \nu_{tur} t}{27} + 
\int_0^t  \|\epsilon(\boldsymbol{\zeta}_{\mathbf{g}}(s)+\mathbf{z}(s)) \|^3_{[L^3(\Omega)]^{3 \times 3}} 
\left\|
\epsilon\left(\frac{\partial \boldsymbol{\zeta}_{\mathbf{g}}}{\partial t}(s)\right)
\right\|_{[L^3(\Omega)]^{3 \times 3}}^{3/2} \, ds  \leq 
\frac{4 \nu_{tur} t}{27} \vspace{0.2cm} \\  \displaystyle + 
\int_0^t \bigg(\|\epsilon(\boldsymbol{\zeta}_{\mathbf{g}}(s)+\mathbf{z}(s))\|_{[L^2(\Omega)]^{3 \times 3}}^2 
+\|\epsilon(\boldsymbol{\zeta}_{\mathbf{g}}(s)+\mathbf{z}(s))\|_{[L^3(\Omega)]^{3 \times 3}}^3 
\bigg) \left\|
\epsilon\left(\frac{\partial \boldsymbol{\zeta}_{\mathbf{g}}}{\partial t}(s)\right)
\right\|_{[L^3(\Omega)]^{3 \times 3}}^{3/2} \, ds.
\end{array} \nonumber
\end{equation}
Taking into account previous estimates:
\begin{equation} \label{eq:inequality1}
\hspace{-1.3cm}\begin{array}{r} 
\displaystyle \frac{1}{2} \left\|\frac{\partial \mathbf{z}}{\partial t} \right\|^2_{L^2(0,t;[L^2(\Omega)]^3)} + 
\nu \|\epsilon(\boldsymbol{\zeta}_{\mathbf{g}}(t)+\mathbf{z}(t))\|^2_{[L^2(\Omega)]^{3 \times 3}}  + 
\frac{2}{3} \nu_{tur} \|\epsilon(\boldsymbol{\zeta}_{\mathbf{g}}(t)+\mathbf{z}(t))\|^3_{[L^3(\Omega)]^{3 \times 3}} 
\vspace{0.2cm} \\  \displaystyle \leq 
\nu \|\epsilon(\mathbf{v}_0) \|^2_{[L^2(\Omega)]^{3 \times 3}} + 
\frac{2}{3} \nu_{tur} \|\epsilon(\mathbf{v}_0) \|^3_{[L^3(\Omega)]^{3 \times 3}} + 
\|\widetilde{\mathbf{H}}_{\mathbf{g}}\|^2_{L^2(0,t;[L^2(\Omega)]^3)}  \vspace{0.2cm} \\  
\displaystyle 
+ C\int_0^t
\bigg( \|\epsilon(\boldsymbol{\zeta}_{\mathbf{g}}(s)+\mathbf{z}(s))\|^2_{[L^{2}(\Omega)]^{3\times 3}}
+\| \epsilon (\boldsymbol{\zeta}_{\mathbf{g}}(s)+\mathbf{z}(s))\|^3_{[L^3(\Omega)]^{3 \times 3}}
\bigg) \| \epsilon (\boldsymbol{\zeta}_{\mathbf{g}}(s)+\mathbf{z}(s))\|^2_{[L^3(\Omega)]^{3 \times 3}} \, ds \vspace{0.2cm} \\ 
\displaystyle+ \frac{\nu t }{4} + \int_0^t 
\bigg(\|\epsilon(\boldsymbol{\zeta}_{\mathbf{g}}(s)+\mathbf{z}(s))\|_{[L^2(\Omega)]^{3 \times 3}}^2 
+\|\epsilon(\boldsymbol{\zeta}_{\mathbf{g}}(s)+\mathbf{z}(s))\|_{[L^3(\Omega)]^{3 \times 3}}^3 
\bigg)
\left\|
\epsilon\left(\frac{\partial \boldsymbol{\zeta}_{\mathbf{g}}}{\partial t}(s)\right)
\right\|_{[L^2(\Omega)]^{3 \times 3}}^2 \, ds
\vspace{0.2cm} \\ \displaystyle +
\frac{4 \nu_{tur} t}{27} + 
\int_0^t \bigg(\|\epsilon(\boldsymbol{\zeta}_{\mathbf{g}}(s)+\mathbf{z}(s))\|_{[L^2(\Omega)]^{3 \times 3}}^2 
+\|\epsilon(\boldsymbol{\zeta}_{\mathbf{g}}(s)+\mathbf{z}(s))\|_{[L^3(\Omega)]^{3 \times 3}}^3 
\bigg) \left\|
\epsilon\left(\frac{\partial \boldsymbol{\zeta}_{\mathbf{g}}}{\partial t}(s)\right)
\right\|_{[L^3(\Omega)]^{3 \times 3}}^{3/2} \, ds.
\end{array} \nonumber
\end{equation}
Then, if we denote by:
\begin{equation} \hspace*{-.9cm}
\begin{array}{rcl}
\Psi_1(s)&=&\displaystyle
\|\epsilon(\boldsymbol{\zeta}_{\mathbf{g}}(s)+\mathbf{z}(s))\|_{[L^2(\Omega)]^{3 \times 3}}^2 
+\|\epsilon(\boldsymbol{\zeta}_{\mathbf{g}}(s)+\mathbf{z}(s))\|_{[L^3(\Omega)]^{3 \times 3}}^3 , \vspace{0.2cm} \\
\Psi_2(s)&=&\displaystyle 
\| \epsilon (\boldsymbol{\zeta}_{\mathbf{g}}(s)+\mathbf{z}(s))\|^2_{[L^3(\Omega)]^{3 \times 3}} +
\left\|
\epsilon\left(\frac{\partial \boldsymbol{\zeta}_{\mathbf{g}}}{\partial t}(s)\right)
\right\|_{[L^2(\Omega)]^{3 \times 3}}^2 
+ \left\|
\epsilon\left(\frac{\partial \boldsymbol{\zeta}_{\mathbf{g}}}{\partial t}(s)\right)
\right\|_{[L^3(\Omega)]^{3 \times 3}}^{3/2} \vspace{0.2cm} \\ & \leq & 
\displaystyle \frac{5}{6} + 
\frac{2}{3} \| \epsilon (\boldsymbol{\zeta}_{\mathbf{g}}(s)+\mathbf{z}(s))\|^3_{[L^3(\Omega)]^{3 \times 3}} +
\left\|
\epsilon\left(\frac{\partial \boldsymbol{\zeta}_{\mathbf{g}}}{\partial t}(s)\right)
\right\|_{[L^2(\Omega)]^{3 \times 3}}^2 
+ \frac{1}{2}\left\|
\epsilon\left(\frac{\partial \boldsymbol{\zeta}_{\mathbf{g}}}{\partial t}(s)\right)
\right\|_{[L^3(\Omega)]^{3 \times 3}}^{3},
\end{array} \nonumber
\end{equation}
we can rewrite the inequality (\ref{eq:inequality1}) in the following terms:
\begin{equation}
\begin{array}{rcl}
\displaystyle
\left\|\frac{\partial \mathbf{z}}{\partial t} \right\|^2_{L^2(0,t;[L^2(\Omega)]^3)} 
+ \Psi_1(t) &\leq&\displaystyle
\nu \|\epsilon(\mathbf{v}_0) \|^2_{[L^2(\Omega)]^{3 \times 3}} + 
\frac{2}{3} \nu_{tur} \|\epsilon(\mathbf{v}_0) \|^3_{[L^3(\Omega)]^{3 \times 3}} 
\vspace{0.2cm} \\ 
& + & \displaystyle 
\|\widetilde{\mathbf{H}}_{\mathbf{g}}\|^2_{L^2(0,t;[L^2(\Omega)]^3)}+
C \int_0^t \Psi_1(s) \Psi_2(s) \, ds.
\end{array}
\end{equation}
In order to derive the estimate for the time derivative, by applying Gronwall's lemma to previous inequality, 
we need first to check that $\Psi_2 \in L^1(0,T)$. So, we have that:
\begin{equation} \hspace*{-1.cm}
\begin{array}{r}
\displaystyle
\int_0^T \Psi_2 (s) \, ds \leq \displaystyle C\Big( 1+
\|\mathbf{v}_0\|_{[L^2(\Omega)]^3}^2
+
\|\boldsymbol{\zeta}_{\mathbf{g}}\|^3_{L^3(0,T;[W^{1,3}(\Omega)]^3)}
\vspace{0.2cm} \\  \displaystyle
+ \|\mathbf{H}_{\mathbf{g}}\|^2_{L^2(0,T;[L^2(\Omega)]^2)}
+\left\| \frac{\partial \boldsymbol{\zeta}_{\mathbf{g}}}{\partial t}\right\|^2_{L^2(0,T;H^1(\Omega))}
+\left\| \frac{\partial \boldsymbol{\zeta}_{\mathbf{g}}}{\partial t}\right\|^3_{L^2(0,T;W^{1,3}(\Omega))} \Big),
\end{array}
\end{equation} 
where we have used estimate (\ref{eq:acotav1}). Then, by Gronwall's lemma, we have a.e. $t \in ]0,T[$:
\begin{equation} \hspace*{-1.cm}
\begin{array}{r}
\displaystyle
\Psi_1(t) \leq C \Big(\nu \|\epsilon(\mathbf{v}_0) \|^2_{[L^2(\Omega)]^{3 \times 3}} + 
\frac{2}{3} \nu_{tur} \|\epsilon(\mathbf{v}_0) \|^3_{[L^3(\Omega)]^{3 \times 3}} 
+
\|\widetilde{\mathbf{H}}_{\mathbf{g}}\|^2_{L^2(0,T;[L^2(\Omega)]^3)} \Big)
\vspace{0.2cm}
\\
\displaystyle
\times
\exp \Big(\|\mathbf{v}_0\|_{[L^2(\Omega)]^3}^2
+
\|\boldsymbol{\zeta}_{\mathbf{g}}\|^3_{L^3(0,T;[W^{1,3}(\Omega)]^3)}
+ \|\mathbf{H}_{\mathbf{g}}\|^2_{L^2(0,T;[L^2(\Omega)]^2)}
\vspace{0.2cm} \\  \displaystyle
+\left\| \frac{\partial \boldsymbol{\zeta}_{\mathbf{g}}}{\partial t}\right\|^2_{L^2(0,T;H^1(\Omega))}
+\left\| \frac{\partial \boldsymbol{\zeta}_{\mathbf{g}}}{\partial t}\right\|^3_{L^2(0,T;W^{1,3}(\Omega))} \Big)
\end{array}
\end{equation}
and then we obtain the existence of a positive constant $C_2$ such that:
\begin{equation} \hspace*{-1.cm}
\begin{array}{r}
\displaystyle
\left\|\frac{\partial \mathbf{z}}{\partial t} \right\|^2_{L^2(0,T;[L^2(\Omega)]^3)} + 
\|\mathbf{z}\|^2_{L^{\infty}(0,T;[W^{1,2}(\Omega)]^3)}+
\|\mathbf{z}\|^3_{L^{\infty}(0,T;[W^{1,3}(\Omega)]^3)} 
\vspace{0.2cm} \\  \displaystyle
\leq 
C_2 \Big(\nu \|\epsilon(\mathbf{v}_0) \|^2_{[L^2(\Omega)]^{3 \times 3}} + 
\frac{2}{3} \nu_{tur} \|\epsilon(\mathbf{v}_0) \|^3_{[L^3(\Omega)]^{3 \times 3}} 
+
\|\widetilde{\mathbf{H}}_{\mathbf{g}}\|^2_{L^2(0,T;[L^2(\Omega)]^3)} \Big)
\vspace{0.2cm} \\  \displaystyle
\times \Bigg[1+
\exp \Big(\|\mathbf{v}_0\|_{[L^2(\Omega)]^3}^2
+
\|\boldsymbol{\zeta}_{\mathbf{g}}\|^3_{L^3(0,T;[W^{1,3}(\Omega)]^3)}
+ \|\mathbf{H}_{\mathbf{g}}\|^2_{L^2(0,T;[L^2(\Omega)]^2)}
\vspace{0.2cm} \\  \displaystyle
+\left\| \frac{\partial \boldsymbol{\zeta}_{\mathbf{g}}}{\partial t}\right\|^2_{L^2(0,T;H^1(\Omega))}
+\left\| \frac{\partial \boldsymbol{\zeta}_{\mathbf{g}}}{\partial t}\right\|^3_{L^2(0,T;W^{1,3}(\Omega))} \Big)
\Bigg],
\end{array}
\end{equation}
which concludes the proof.
\hfill $\blacksquare$

\subsection{Galerkin approximation}

In this subsection we will construct a sequence of approximations to the solution 
aimed at converging to a solution of problem (\ref{eq:system1}). In order to 
pass to the limit in below approximations, we will need the following technical result:

\begin{lem} \label{mono} The operator:
\begin{equation}
\begin{array}{rcl}
A:\widetilde{\mathbf{X}}& \rightarrow &\widetilde{\mathbf{X}}' \\ 
\mathbf{z} & \rightarrow & A(\mathbf{z}),
\end{array}
\end{equation}
where, for any $\boldsymbol{\xi}\in \widetilde{\mathbf{X}}$,
\begin{equation}
\langle A(\mathbf{z}),\boldsymbol{\xi} \rangle =\int_{\Omega} 
\beta(\epsilon(\boldsymbol{\zeta}_{\mathbf{g}}+\mathbf{z})) 
\epsilon(\boldsymbol{\zeta}_{\mathbf{g}}+\mathbf{z}):\epsilon(\boldsymbol{\xi}) \, 
d \mathbf{x},
\end{equation}
satisfies the following monotony condition:
\begin{equation}
\langle A(\mathbf{z}_1)-A(\mathbf{z}_2),\mathbf{z}_1-\mathbf{z}_2 \rangle \geq 
C \int_{\Omega} \|\nabla (\mathbf{z}_1-\mathbf{z}_2)\|_{[L^2(\Omega)]^{3 \times 3}}^2. 
\end{equation}
\end{lem}

\noindent {\bf Proof } Using classical results of integral calculus:
\begin{equation} \label{eq:monotono1}
\begin{array}{c}
\beta (\epsilon(\mathbf{v}_1)) \epsilon (\mathbf{v}_1) 
-\beta (\epsilon(\mathbf{v}_2)) \epsilon (\mathbf{v}_2) =
\displaystyle 
\left. \frac{\partial D(\epsilon)}{\partial \epsilon} \right|_{\epsilon=\epsilon(\mathbf{v}_1)}-
\left. \frac{\partial D(\epsilon)}{\partial \epsilon} \right|_{\epsilon=\epsilon(\mathbf{v}_2)} \\
= \displaystyle 
\int_0^1 \frac{d}{d \tau} \left. \frac{\partial D(\epsilon)}{\partial \epsilon} 
\right|_{\epsilon=\epsilon_{\tau}}\, d \tau 
= \displaystyle \int_0^1 \left. \frac{\partial^2 D(\epsilon)}{\partial \epsilon^2} 
\right|_{\epsilon=\epsilon_{\tau}} \, d \tau \left(\epsilon(\mathbf{v}) \right),
\end{array}
\end{equation}
where $\mathbf{v}_1$ and $\mathbf{v}_2$ are two elements of $\mathbf{X}$, 
$\epsilon_{\tau}=\tau \epsilon(\mathbf{v}_1)+(1-\tau) \epsilon(\mathbf{v}_2)$, with 
$\tau \in [0,1]$,and, consequently, $\frac{\partial}{\partial \tau} \epsilon_{\tau} = 
\epsilon(\mathbf{v}_1)-\epsilon(\mathbf{v}_2) =\epsilon(\mathbf{v})$. Then, taking
into account the definition (\ref{eq:doperador}) of operator $D$, it is straightforward 
to check that:
\begin{equation} \label{eq:monotono2} 
\begin{array}{r}
\displaystyle 
\left(
\int_0^1 \left. \frac{\partial^2 D(\epsilon)}{\partial \epsilon^2} 
\right|_{\epsilon=\epsilon_{\tau}} \, d \tau 
\left(\epsilon(\mathbf{v}) \right),\epsilon(\mathbf{v})\right)_{[L^2(\Omega)]^{3 \times 3}}
\\ \displaystyle = 
\int_0^1 \left(
 \left. \frac{\partial^2 D(\epsilon)}{\partial \epsilon^2} 
\right|_{\epsilon=\epsilon_{\tau}} 
\left(\epsilon(\mathbf{v}) \right),\epsilon(\mathbf{v})\right)_{[L^2(\Omega)]^{3 \times 3}}
\, d \tau \\ \displaystyle
 \geq 
 C_1\|\epsilon(\mathbf{v})\|^2_{[L^2(\Omega)]^{3 \times 3}}
  \int_0^1 \left(1+
\|\epsilon_{\tau}\|_{[L^1(\Omega)]^{3 \times 3}}
\right)\, d \tau  
\geq  
C_2  \|\epsilon(\mathbf{v})\|^2_{[L^2(\Omega)]^{3 \times 3}},
\end{array}
\end{equation}
where $C_1$ y $C_2$ are positive constants. Thus, we can conclude that there 
exists a positive constant $C$ such that:
\begin{equation}
\langle A(\mathbf{z}_1)-A(\mathbf{z}_2),\mathbf{z}_1-\mathbf{z}_2 \rangle \geq 
C \int_{\Omega} \|\nabla (\mathbf{z}_1-\mathbf{z}_2)\|_{[L^2(\Omega)]^{3 \times 3}}^2.
\end{equation}
\hfill $\blacksquare$

\begin{thm} There exists a solution $\mathbf{z} \in \widetilde{\mathbf{W}}$ 
of the system (\ref{eq:system2}) in the sense of Definition \ref{defsol}.
\end{thm}

\noindent {\bf Proof } We consider a basis $\{\boldsymbol{\xi}^n\}_{n \in \mathbb{N}}$ 
of the functional space $\{\mathbf{v} \in W^{2,2}(\Omega):\; \nabla \cdot \mathbf{v}=0\; \mbox{and}\; 
\mathbf{v}_{|_{\partial \Omega}}=\boldsymbol{0}\}$, formed by eigenfunctions of the Stokes operator:
\begin{equation} 
\left\{\begin{array}{l}
\displaystyle
-\Delta \boldsymbol{\xi}^n + \nabla p^n= \lambda^n \boldsymbol{\xi}^n \quad \mbox{in}\; \Omega, 
\\
\displaystyle
\nabla \cdot \boldsymbol{\xi}^n=0 \quad \mbox{in} \; \Omega,
\\
\displaystyle
\boldsymbol{\xi}^n=0 \quad \mbox{on}\; \partial \Omega,
\end{array}\right.
\end{equation}
where $0<\lambda_1\leq \lambda_2 \leq \ldots $ are the eigenvalues, with $\lim_{n \to \infty} \lambda_n=\infty$. We 
can also suppose that the previous basis is orthonormal with respect to $L^2(\Omega)$. For $N\in \mathbb{N}$ we denote 
by:
\begin{equation}
\mathbf{z}^N=\sum_{k=1}^N \mathbf{z}_k^N(t) \boldsymbol{\xi}^k,
\end{equation}
where the coefficients $\mathbf{z}_k^N(t)$, $k=1,\ldots,N$, are such that $\mathbf{z}^N$ is the solution of the 
following differential equation, for $k=1,\ldots,N$:
\begin{equation} \label{eq:system4}
\hspace{-1.1cm} \begin{array}{r}
\displaystyle
\int_{\Omega} \frac{\partial \mathbf{z}^N}{\partial t} \cdot \boldsymbol{\xi}^k \,d\mathbf{x} +
\int_{\Omega} \nabla (\boldsymbol{\zeta}_{\mathbf{g}}+
\mathbf{z}^N) \mathbf{z}^N \cdot \boldsymbol{\xi}^k \, d \mathbf{x} \vspace{0.2cm} \\
+ \displaystyle
\int_{\Omega} \nabla \mathbf{z}^N \boldsymbol{\zeta}_{\mathbf{g}} 
\cdot \boldsymbol{\xi}^k \, d \mathbf{x}
+ 2 \nu \int_{\Omega} \epsilon(\mathbf{z}^N) : \epsilon(\boldsymbol{\xi}^k) \, 
d \mathbf{x} \vspace{0.2cm} \\
+ \displaystyle
2 \nu_{tur} \int_{\Omega} \left[\epsilon(\boldsymbol{\zeta}_{\mathbf{g}}+\mathbf{z}^N):
\epsilon(\boldsymbol{\zeta}_{\mathbf{g}}+\mathbf{z}^N) \right]^{1/2} \epsilon(\boldsymbol{\zeta}_{\mathbf{g}}+\mathbf{z}^N):
\epsilon(\boldsymbol{\xi}^k) \, d \mathbf{x} 
=\int_{\Omega} \mathbf{H}_{\mathbf{g}} \cdot \boldsymbol{\xi}^k \, d \mathbf{x},
\end{array}
\end{equation}
which can be rewritten in the following standard matrix formulation:
\begin{equation} \label{eq:system5}
\left\{\begin{array}{l}
\displaystyle \frac{d\mathbf{z}(t)}{dt}=\mathbf{F}(\mathbf{z}(t),t) \quad \mbox{a.e. } \, t \in ]0,T[,  \\ 
\displaystyle \mathbf{z}(0)=\mathbf{z}_0,
\end{array}\right.
\end{equation}
where:
\begin{equation}
\mathbf{z}(t)=\left(\mathbf{z}_1^N(t),\, \mathbf{z}_2^N(t),\ldots,\mathbf{z}_N^N(t) \right)^T,
\end{equation}
\begin{equation}
\mathbf{z}_0=\left(
(\mathbf{z}_0^N,\boldsymbol{\xi}^1),
(\mathbf{z}_0^N,\boldsymbol{\xi}^2),\ldots, 
(\mathbf{z}_0^N,\boldsymbol{\xi}^N) \right)^T,
\end{equation}
\begin{equation}
\mathbf{F}(\mathbf{z},t)=\left(\begin{array}{c}
(\mathbf{H}_{\mathbf{g}}(t),\boldsymbol{\xi}^1)-a(t;\mathbf{z} \cdot \boldsymbol{\xi},\boldsymbol{\xi}^1) \\
(\mathbf{H}_{\mathbf{g}}(t),\boldsymbol{\xi}^2)-a(t;\mathbf{z} \cdot \boldsymbol{\xi},\boldsymbol{\xi}^2) \\
\vdots \\
(\mathbf{H}_{\mathbf{g}}(t),\boldsymbol{\xi}^N)-a(t;\mathbf{z} \cdot \boldsymbol{\xi},\boldsymbol{\xi}^N)
\end{array}\right),
\end{equation}
\begin{equation}
\boldsymbol{\xi}=\left( 
\boldsymbol{\xi}^1,\boldsymbol{\xi}^2,\ldots,\boldsymbol{\xi}^N
\right)^T,
\end{equation}
\begin{equation}\hspace*{-.9cm}
\begin{array}{r}
\displaystyle a(t;\mathbf{z} \cdot \boldsymbol{\xi},\boldsymbol{\xi}^k) =
\displaystyle 
\int_{\Omega} \nabla (\boldsymbol{\zeta}_{\mathbf{g}}(t)+
(\mathbf{z} \cdot \boldsymbol{\xi})) (\mathbf{z} \cdot \boldsymbol{\xi}) 
\cdot \boldsymbol{\xi}^k \, d \mathbf{x} +
\int_{\Omega} \nabla (\mathbf{z} \cdot \boldsymbol{\xi}) \boldsymbol{\zeta}_{\mathbf{g}} 
\cdot \boldsymbol{\xi}^k \, d \mathbf{x} \vspace{0.2cm}\\  \displaystyle+
2 \nu_{tur} \int_{\Omega} \left[\epsilon(\boldsymbol{\zeta}_{\mathbf{g}}+(\mathbf{z} \cdot \boldsymbol{\xi})):
\epsilon(\boldsymbol{\zeta}_{\mathbf{g}}(t)+(\mathbf{z} \cdot \boldsymbol{\xi})) \right]^{1/2} 
\epsilon(\boldsymbol{\zeta}_{\mathbf{g}}(t)+(\mathbf{z} \cdot \boldsymbol{\xi})):
\epsilon(\boldsymbol{\xi}^k) \, d \mathbf{x}
 \vspace{0.2cm}\\ 
 \displaystyle
+ 2 \nu \int_{\Omega} \epsilon(\mathbf{z} \cdot \boldsymbol{\xi}) : \epsilon(\boldsymbol{\xi}^k) \, 
d \mathbf{x},\quad k=1,\ldots,N,
\end{array}
\end{equation}
for $\mathbf{z}_0^N$ the orthogonal projection in $L^2(\Omega)$ of 
$\mathbf{v}_0$ onto $\langle \{\boldsymbol{\xi}_1,\ldots, 
\boldsymbol{\xi}_N\} \rangle$. 

Then, we can apply the Caratheodory theorem to this system of ordinary differential equations
(cf., for instance, Theorem 5.2 of \cite{hale1}). Indeed, $\mathbf{F}(\cdot,t)$ is 
continuous for any $t \in [0,T]$, and $\mathbf{F}(\mathbf{z}, \cdot ) \in 
L^{2}(0,T)$ for any $\mathbf{z} \in \mathbb{R}^{3N}$. Thus, given an open 
ball $B$ in $\mathbb{R}^{3N}$, we can prove that there exist two functions 
$m_B,\, l_B \in L^1(0,T)$ (in fact, in $L^2(0,T)$) such that:
\begin{equation} \label{eq:caratheocond}
\begin{array}{rcl}
\displaystyle 
\|\mathbf{F}(\mathbf{z},t)\| &\leq & \displaystyle 
m_B(t) \quad \mbox{a.e.} \; t \in ]0,T[, \quad \forall\, \mathbf{z} \in B, \\
\displaystyle
\|\mathbf{F}(\mathbf{z}_1,t)-
\mathbf{F}(\mathbf{z}_2,t)\| & \leq & 
\displaystyle l_B(t) \|\mathbf{z}_1-\mathbf{z}_2\| \quad \mbox{a.e.} \; t \in ]0,T[, \quad
\forall\, \mathbf{z}_1,\, \mathbf{z}_2 \in B,
\end{array}
\end{equation}
and then, we conclude that system (\ref{eq:system5}) has a unique 
absolutely continuous solution which can be extended to the boundary 
of $B \times ]0,T[ $. Moreover, thanks to the 
regularity of $\mathbf{F}$ in time, $\mathbf{z}_k^N \in H^1(0,T)$, $\forall k=1,
\ldots,N$, and then, $\mathbf{z}^N \in W^{1,2,2}(0,T;[H^2(\Omega)],[H^2(\Omega)]^3)$. 

It is important to mention here that we can repeat the 
proof of the estimates (\ref{eq:acotav1}) and (\ref{eq:acotav2}) in 
the Galerkin approximation (\ref{eq:system4}) and then the sequence 
$\{\mathbf{z}^N\}_{N \in \mathbb{N}}$ is bounded in the space 
$\widetilde{\mathbf{W}}$ by a constant not depending on $N$. 

The final part of this proof corresponds to the pass to the limit in the 
Galerkin approximations in order to obtain a solution of the system 
(\ref{eq:system2}). Taking subsequences, if necessary, we have that:
\begin{itemize}
\item $\mathbf{z}^N \rightarrow \mathbf{z} $ strongly in 
$L^p(0,T;[L^q(\Omega)]^3)$, for all $1<p<\infty$ and $2\leq q < \infty$,
\item $\displaystyle \mathbf{z}^N \rightharpoonup \mathbf{z}$ weakly in 
$L^3(0,T;\widetilde{\mathbf{X}})$,
\item $\displaystyle \frac{d \mathbf{z}^N}{dt} \rightharpoonup 
\frac{d \mathbf{z}}{dt}$ weakly in $L^2(0,T;[L^2(\Omega)]^3)$,
\item $\displaystyle \nabla \mathbf{z}^N \rightharpoonup^* \nabla 
\mathbf{z}$ weakly-$*$ in $L^{\infty}(0,T;[L^3(\Omega)]^3)$,
\item $\displaystyle \beta(\epsilon(\boldsymbol{\zeta}_{\mathbf{g}}+\mathbf{z}^N)) \, 
\epsilon(\boldsymbol{\zeta}_{\mathbf{g}}+
\mathbf{z}^N)\rightharpoonup \widehat{\beta}$ weakly 
in $L^{3/2}(0,T;\widetilde{\mathbf{X}}')$.
\end{itemize}

Let us fix now an index $k \in \mathbb{N}$. If we multiply (\ref{eq:system4}) by a scalar 
function $\psi\in H^1(0,T)$, such that $\psi(T)=0$, 
integrate with respect to $t$, and integrate by parts, we have, $\forall N \geq k$:
\begin{equation} \label{eq:system7} \hspace*{-.8cm}
\begin{array}{r}
\displaystyle -\int_0^T \int_{\Omega} \mathbf{z}^N(t) \cdot
\frac{d\psi}{dt}(t)  \boldsymbol{\xi}^k\, d \mathbf{x}\, dt 
+\int_0^T \int_{\Omega} \nabla (\boldsymbol{\zeta}_{\mathbf{g}}+\mathbf{z}^N)
(\boldsymbol{\zeta}_{\mathbf{g}}+\mathbf{z}^N) \cdot \psi(t) \boldsymbol{\xi}^k \, 
d \mathbf{x}\, dt\vspace{0.2cm} \\ 
\displaystyle +\int_0^T \int_{\Omega} 
\beta(\epsilon(\boldsymbol{\zeta}_{\mathbf{g}}+\mathbf{z}^N))
\epsilon(\boldsymbol{\zeta}_{\mathbf{g}}+\mathbf{z}^N):\epsilon(\psi(t)
\boldsymbol{\xi}^k) \, d \mathbf{x}\, dt
=\int_{\Omega} \mathbf{z}^N_0 \cdot \psi(0) \boldsymbol{\xi}^k \, d \mathbf{x}
\vspace{0.2cm} \\  \displaystyle 
+\int_0^T \int_{\Omega} 
\mathbf{H}_{\mathbf{g}} \cdot \psi(t) \boldsymbol{\xi}^k \, d\mathbf{x} \, dt.
\end{array}
\end{equation}
Thanks to the previous convergences we can pass to the limit in this 
expression and we obtain:
\begin{equation}\label{eq:system6}
\begin{array}{rcl}
\displaystyle -\int_0^T \int_{\Omega} \mathbf{z}(t) \cdot
\frac{d\psi}{dt}(t)  \boldsymbol{\xi}\, d \mathbf{x}\, dt 
+\int_0^T \int_{\Omega} \nabla (\boldsymbol{\zeta}_{\mathbf{g}}+\mathbf{z})
(\boldsymbol{\zeta}_{\mathbf{g}}+\mathbf{z}) \cdot \psi(t) \boldsymbol{\xi} \, 
d \mathbf{x}\, dt\vspace{0.2cm} \\ 
\displaystyle +\int_0^T \int_{\Omega} 
\widehat{\beta}:\epsilon(\psi(t)
\boldsymbol{\xi}) \, d \mathbf{x}\, dt
= \int_{\Omega} \mathbf{z}_0 \cdot \psi(0) \boldsymbol{\xi} \, d \mathbf{x}
+ \int_0^T \int_{\Omega} 
\mathbf{H}_{\mathbf{g}} \cdot \psi(t) \boldsymbol{\xi}\, d\mathbf{x} \, dt,
\end{array}
\end{equation}
for each $ \boldsymbol{\xi}\in \widetilde{\mathbf{X}}$ which is a finite linear 
combination of elements $\boldsymbol{\xi}^k$. Since each term 
of above expression depends linearly and continuously on 
$\boldsymbol{\xi}$ for the norm of $\widetilde{\mathbf{X}}$, 
previous equality is still valid, by continuity, for each $\boldsymbol{\xi} 
\in \widetilde{\mathbf{X}}$. Then, we are only left to demonstrate that 
$\widehat{\beta} = \beta(\epsilon(\boldsymbol{\zeta}_{\mathbf{g}}+\mathbf{z})) \, 
\epsilon(\boldsymbol{\zeta}_{\mathbf{g}}+
\mathbf{z})$. 

So, by Lemma \ref{mono}, we have:
\begin{equation}
\begin{array}{rcl}
\displaystyle 
\int_0^T \int_{\Omega}\left[
 \beta(\epsilon(\boldsymbol{\zeta}_{\mathbf{g}}+\mathbf{z}^N)) \, 
\epsilon(\boldsymbol{\zeta}_{\mathbf{g}}+
\mathbf{z}^N)-
 \beta(\epsilon(\boldsymbol{\zeta}_{\mathbf{g}}+
 \boldsymbol{\eta})) 
\epsilon(\boldsymbol{\zeta}_{\mathbf{g}}+
\boldsymbol{\eta})\right] :  \epsilon(\mathbf{z}^N-
\boldsymbol{\eta}) \, d \mathbf{x}\, dt \geq 0,
\end{array} \nonumber
\end{equation}
for all $\boldsymbol{\eta}(t,\mathbf{x})=
\psi(t) \boldsymbol{\xi}(\mathbf{x})$, with $\psi\in H^1(0,T)$ and $\boldsymbol{\xi} \in 
\langle\{\boldsymbol{\xi}_1,\ldots,\boldsymbol{\xi}_N\}\rangle$. Then, by (\ref{eq:system7}), 
\begin{equation} \hspace*{-1.4cm}
\begin{array}{r}
\displaystyle 
\int_0^T \int_{\Omega} \mathbf{z}^N(t) \cdot \frac{d}{dt} (\mathbf{z}^N-
\boldsymbol{\eta})\,d \mathbf{x}\, dt -
\int_0^T \int_{\Omega} \nabla (\boldsymbol{\zeta}_{\mathbf{g}}+\mathbf{z}^N)
(\boldsymbol{\zeta}_{\mathbf{g}}+\mathbf{z}^N) \cdot
 (\mathbf{z}^N-
\boldsymbol{\eta})\,d \mathbf{x}\, dt
\vspace{0.2cm} \\ \displaystyle
+\int_0^T  \mathbf{z}_0^N \cdot  (\mathbf{z}^N(0)-
\boldsymbol{\eta}(0))\,d \mathbf{x}\, dt
+\int_0^T \int_{\Omega} \mathbf{H}_{\mathbf{g}} \cdot 
 (\mathbf{z}^N-
\boldsymbol{\eta})\,d \mathbf{x}\, dt
\vspace{0.2cm} \\ \displaystyle
-\int_0^T \int_{\Omega}
 \beta(\epsilon(\boldsymbol{\zeta}_{\mathbf{g}}+
 \boldsymbol{\eta})) 
\epsilon(\boldsymbol{\zeta}_{\mathbf{g}}+
\boldsymbol{\eta}) :  \epsilon(\mathbf{z}^N-
\boldsymbol{\eta}) \, d \mathbf{x}\, dt  \geq 0.
\end{array}
\end{equation}
If we pass to the limit, taking into account (\ref{eq:system6}):
\begin{equation}
\int_0^T \int_{\Omega} \left[ \widehat{\beta}-
 \beta(\epsilon(\boldsymbol{\zeta}_{\mathbf{g}}+
 \boldsymbol{\eta})) 
\epsilon(\boldsymbol{\zeta}_{\mathbf{g}}+
\boldsymbol{\eta})\right] :  \epsilon(\mathbf{z}-
\boldsymbol{\eta})  \, d \mathbf{x}\, dt  \geq 0,
\end{equation}
for all $\boldsymbol{\eta}\in L^3(0,T;\widetilde{\mathbf{X}})$ 
(where we need to use a standard argument by density). Then, choosing 
$\boldsymbol{\eta}=\mathbf{z} \pm \epsilon \boldsymbol{\zeta}$, with 
$\boldsymbol{\zeta} \in  L^3(0,T;\widetilde{\mathbf{X}})$, and   
$\epsilon$ an arbitrary positive number, multiplying both sides of 
the inequality by $\epsilon^{-1}$, and letting $\epsilon$ tend to zero, 
we obtain that, for all $\boldsymbol{\zeta}$:
\begin{equation}
\int_0^T \int_{\Omega} 
\left[ \widehat{\beta}-
 \beta(\epsilon(\boldsymbol{\zeta}_{\mathbf{g}}+
 \mathbf{z})) 
\epsilon(\boldsymbol{\zeta}_{\mathbf{g}}+
\mathbf{z})\right] :  
\boldsymbol{\zeta}  \, d \mathbf{x}\, dt  =0.
\end{equation}
Thus, $\widehat{\beta}=
 \beta(\epsilon(\boldsymbol{\zeta}_{\mathbf{g}}+
 \mathbf{z})) 
\epsilon(\boldsymbol{\zeta}_{\mathbf{g}}+
\mathbf{z})$ a.e. $(\mathbf{x},t) \in \Omega \times ]0,T[$. 

Finally, since we can derive directly from the definition of $\mathbf{z}_0$ that 
$\mathbf{z}(0)=\mathbf{v}_0$, the proof is complete.
\hfill $\blacksquare$

\section{Uniqueness of solution}

\begin{thm} Under the hypotheses that guaranty the existence of solution for the problem (\ref{eq:system2}), 
there exists a unique solution in the sense of Definition \ref{defsol}. 
\end{thm}

\noindent {\bf Proof } 
We suppose that there exist $\mathbf{v}_1$ 
and $\mathbf{v}_2$ two solutions of (\ref{eq:system2}) in the sense of 
definition \ref{defsol}, that is, there exists $\mathbf{z}_1$ and $\mathbf{z}_2$ 
in the space $\widetilde{\mathbf{W}}$ such that 
$\mathbf{v}_1=\boldsymbol{\zeta}_{\mathbf{g}}+
\mathbf{z}_1$ and $\mathbf{v}_2=\boldsymbol{\zeta}_{\mathbf{g}}+
\mathbf{z}_2$ and satisfy the variational formulation (\ref{eq:system3}). 
If we denote by $\mathbf{v}=\mathbf{z}_1-\mathbf{z}_2=\mathbf{v}_1
-\mathbf{v}_2\in \widetilde{\mathbf{W}}$, we have that $\mathbf{v}(0)=\boldsymbol{0}$, and:
\begin{equation}\hspace*{-.9cm}
\begin{array}{r}
\displaystyle 
\int_{\Omega} \frac{\partial \mathbf{v}}{\partial t} \cdot 
\boldsymbol{\eta} \,d\mathbf{x} +
\int_{\Omega} \left( 
\nabla \mathbf{v}_1 \mathbf{v}_1 
-\nabla \mathbf{v}_2 \mathbf{v}_2
\right) \cdot \boldsymbol{\eta} \, d \mathbf{x} \\ 
\displaystyle 
+ \int_{\Omega} \left(\beta (\epsilon(\mathbf{v}_1)) \epsilon (\mathbf{v}_1) 
-\beta (\epsilon(\mathbf{v}_2)) \epsilon (\mathbf{v}_2) \right)
: \epsilon(\boldsymbol{\eta})\, d \mathbf{x} =0,\;
\mbox{a.e.} \ t \in ]0,T[, \ \forall
 \boldsymbol{\eta} \in \widetilde{\mathbf{X}}.
\end{array}
\end{equation}
On one hand, we know that
\begin{equation}
\nabla \mathbf{v}_1 \mathbf{v}_1 
-\nabla \mathbf{v}_2 \mathbf{v}_2 =
\nabla \mathbf{v}_1 \mathbf{v}+ \nabla \mathbf{v} \mathbf{v}_2
\end{equation}
and, if we take $\boldsymbol{\eta}=\mathbf{v}(t)$ as a test function and use 
that $\nabla \cdot \mathbf{v}=0$, we obtain:
\begin{equation}
\begin{array}{rcl}
\displaystyle \left( \nabla \mathbf{v} \mathbf{v}_2, \mathbf{v} \right)_{[L^2(\Omega)]^3} &=&0, \\
\displaystyle \left( \nabla \mathbf{v}_1 \mathbf{v},\mathbf{v} \right)_{[L^2(\Omega)]^3} &=&
\displaystyle -\left(\nabla \mathbf{v} \mathbf{v} , \mathbf{v}_1 \right)_{[L^2(\Omega)]^3}.
\end{array}
\end{equation}
Then, thanks to Lemma \ref{mono}:
\begin{equation}
\begin{array}{r}
\displaystyle 
\frac{d}{dt} \| \mathbf{v}(t)\|^2_{[L^2(\Omega)]^3}+
2C 
\|\nabla(\mathbf{v})(t)\|^2_{[L^2(\Omega)]^{3 \times 3}}
\leq 
 \left(\nabla \mathbf{v}(t) \mathbf{v}(t) , \mathbf{v}_1(t) \right)_{[L^2(\Omega)]^3}.
\end{array}
\end{equation}
Using now the fact that $\mathbf{W} \subset L^p(0,T;[L^q(\Omega)]^3)$ for all 
$1<p,\,q<\infty$:
\begin{equation}
\begin{array}{r}
\displaystyle
 \left(\nabla \mathbf{v}(t) \mathbf{v}(t) , \mathbf{v}_1(t) \right)_{[L^2(\Omega)]^3}
 \leq \displaystyle
 \|\nabla \mathbf{v}(t)\|_{[L^2(\Omega)]^{3 \times 3}} 
 \|\mathbf{v}(t)\|_{[L^4(\Omega)]^3}
 \|\mathbf{v}_1(t)\|_{[L^4(\Omega)]^3} \vspace{0.2cm} \\ 
 \leq  C_1
 \displaystyle 
  \|\nabla \mathbf{v}(t)\|_{[L^2(\Omega)]^{3 \times 3}}^{7/4}
 \|\mathbf{v}(t)\|_{[L^2(\Omega)]^3}^{1/4}
 \|\mathbf{v}_1(t)\|_{[L^4(\Omega)]^3} 
  \vspace{0.2cm} \\ 
 \leq  \displaystyle 2 C 
 \|\nabla(\mathbf{v})(t)\|^2_{[L^2(\Omega)]^{3 \times 3}} +
 C_2  \|\mathbf{v}(t)\|_{[L^2(\Omega)]^3}^2
  \|\mathbf{v}_1(t)\|_{[L^4(\Omega)]^3}^8,
\end{array}
\end{equation}
where the second inequality is a consequence of lemma 3.5 of \cite{Temam} and 
the third one is a direct consequence of Young inequality. Finally, we deduce that:
\begin{equation}
\frac{d}{dt} \| \mathbf{v}(t)\|^2_{[L^2(\Omega)]^3} \leq 
C_2  \|\mathbf{v}(t)\|_{[L^2(\Omega)]^3}^2
  \|\mathbf{v}_1(t)\|_{[L^4(\Omega)]^3}^8
\end{equation}
and, since the function $t \rightarrow  \|\mathbf{v}_1(t)\|_{[L^4(\Omega)]^3}^8$ 
is integrable, we can multiply the previous inequality by 
$\exp \left(
-C_2 \int_0^t 
\|\mathbf{v}_1(t)\|_{[L^4(\Omega)]^3}^8\, dt
\right)$, and we obtain:
\begin{equation}
\frac{d}{dt} \left[
\exp \left(
-C_2 \int_0^t 
\|\mathbf{v}_1(t)\|_{[L^4(\Omega)]^3}^8\, dt
\right)\| \mathbf{v}(t)\|^2_{[L^2(\Omega)]^3}
\right] \leq 0.
\end{equation}
Finally, integrating and taking into account that $\mathbf{v}(0)=\boldsymbol{0}$, we find that:
\begin{equation}
\|\mathbf{v}(t)\|^2_{[L^2(\Omega)]^3} \leq 0, \quad \forall t \in [0,T],
\end{equation}
which implies that $\mathbf{v}=\boldsymbol{0}$ or, equivalently, that $\mathbf{v}_1=\mathbf{v}_2$.
\hfill $\blacksquare$



\end{document}